\documentclass[a4paper]{article}

\usepackage{MyLaTeX,Myarticle,graphicx}
\title{Trivial-source endotrivial modules for sporadic groups}
\author{David A. Craven}
\begin{document}
\maketitle

\begin{abstract} We determine the group of endotrivial modules (as an abstract group) for $G$ a (quasi)simple group of sporadic type, extending previous results in the literature. In many sporadic cases we directly construct the subgroup of trivial-source endotrivial modules. We also resolve the question of whether certain simple modules for sporadic groups are endotrivial, posed by Lassueur, Malle and Schulte, in the majority of open cases. The results rely heavily on a recent description of the group of trivial-source endotrivial modules due to Grodal.
\end{abstract}

\section{Introduction}

Dade once wrote ``There are just too many modules over $p$-groups'', and introduced endopermutation -- modules $M$ such that $M\otimes M^*$ is a permutation module -- and endotrivial  -- where it is the sum of a projective and the trivial module -- modules as a consequence \cite{dade1978}. These modules have found their way into many parts of representation theory, from the structure of nilpotent blocks to equivalences of categories. It is fair to say that the endotrivial modules are one of the most important classes of modules in the representation theory of finite groups.

The set of isomorphism classes of indecomposable endotrivial modules for a finite group $G$ possesses a natural abelian group structure induced by the tensor product, called the group of endotrivial modules, and denoted $T(G)$. If $G$ is a finite $p$-group then the structure and generators for $T(G)$ have been completely determined \cite{carlsonthevenaz2004}. Over the last decade, the focus has shifted to understanding $T(G)$ for $G$ an arbitrary finite group. Any attempt to apply the classification of finite simple groups and Clifford theory will necessarily require a detailed understanding of $T(G)$ for $G$ a finite simple group with Sylow $p$-subgroup $S$. The restriction of an endotrivial module is always endotrivial, and so there is a group homomorphism $T(G)\to T(S)$, which by \cite{carlsonmazzathevenaz2013,mazzathevenaz2007} is a split map. Write $T(G,S)$ for the kernel, the subgroup of $T(G)$ consisting of trivial-source endotrivial modules. The determination of $T(G,S)$ is the main problem in endotrivial modules for arbitrary finite groups.

In \cite{lassueurmazza2015}, Lassueur and Mazza made some headway on the problem of determining $T(G,S)$ for $G$ a quasisimple sporadic group by computing it for small sporadic groups (the Monster was considered in \cite{grodal2016un}), and in \cite{carlsonhemmermazza2010,carlsonmazzanakano2009} Carlson--Hemmer--Mazza--Nakano determine $T(G,S)$ for $G$ an alternating or symmetric group (their covering groups were studied in \cite{lassueurmazza2015a}). Both results contain minor errors ($T(G,S)$ for $G=3.J_3$ and $p=2$ should be $1$, not $3$, and for $G=A_{2p},A_{2p+1}$ and $p>3$ it should be $4$, not $2\times 2$), the former corrected here and the latter in \cite{grodal2016un}. We complete the determination of $T(G,S)$ for all sporadic groups $G$, and determine $T(G,S)$ for $G$ `almost quasisimple' of sporadic type, e.g., a double cover of a sporadic. Together with work recently announced by Carlson, Grodal, Mazza and Nakano (and still in preparation) that completes the determination of $T(G,S)$ for Lie type groups in non-defining characteristic, this completes the determination of $T(G,S)$ for all finite quasisimple groups $G$.

We show that $T(G,S)\neq 1$ if $G=J_3$ and $p=3$. This group has $3$-rank $3$, and is the only example here of a quasisimple group $G$ having $p$-rank greater than $2$ and with $T(G,S)\neq 1$. This compares with \cite{lassueurmalle2015}, where it is shown that if a quasisimple group $G$ possesses a simple endotrivial module then it has $p$-rank at most $2$. (The groups $G=\SL_2(p^a)$ have $T(G,S)\cong C_{p^a-1}$, because they have Sylow $p$-subgroups that intersect trivially (see \cite{carlsonmazzanakano2006}), and of course have arbitrarily large $p$-rank.)

\begin{thm}\label{thm:mainthm} Let $G$ be a quasisimple sporadic group or the Tits group $T={}^2\!F_4(2)'$. The group $T(G,S)$ is the abelianization (quotient by the derived subgroup) of the group in Table \ref{tab:result}.

If $G$ is a quasisimple sporadic group such that the $p$-rank of $G$ is at least $3$ and $T(G,S)\neq 1$, then $G=J_3$ and $p=3$ (which has $p$-rank exactly $3$).
\end{thm}
\begin{table}
\begin{center}
\begin{tabular}{ccccccc}
\hline Group & $p=2$ & $p=3$ & $p=5$ & $p=7$ & $p=11$ & $p=13$
\\\hline  $M_{11}$ & $1$ (SN) & $SD_{16}$ (TI) &&&&
\\ $M_{12}$ & $1$ (SN) & $1$ ($K^\circ$) &&&&
\\ $M_{22}$ & $1$ (SN) & $Q_8$ (r2) &&&&
\\ $2\cdot M_{22}$ & $-$ & SmallGroup(16,4)&&&&
\\ $M_{23}$ & $1$ (SN) & $2$ (r2) &&&&
\\ $M_{24}$ & $1$ (SN) & $1$ ($K^\circ$) &&&&
\\ $J_1$ & $1$ ($K^\circ$) &&&&&
\\ $J_2$ & $1$ ($K^\circ$) & $2$ (nnc) & $2$ (r2) &&&
\\ $J_3$ & $1$ ($K^\circ$) & $2$ (nnc)$^*$ &&&&
\\ $J_4$ & $1$ (SN) & $1$ ($K^\circ$) &&&$5\times 2S_4$ (TI)&
\\ $HS$ & $1$ (SN) & $2\times 2$ (r2) & $4$ (nnc) &&&
\\ $2\cdot HS$ & $-$ & $D_8$ & $-$ &&&
\\ $McL$ & $1$ (SN) & $1$ ($K^\circ$) & $3\rtimes 8$ (TI) &&&
\\ $3\cdot McL$ & $-$ & $-$ & $3\times (3\rtimes 8)$&&&
\\ $Suz$ & $1$ ($K^\circ$) & $1$ ($K^\circ$) & $2$ (r2) &&&
\\ $He$ & $1$ (SN) & $1$ ($K^\circ$) & $3$ (r2) & $1$ ($K^\circ$)$^*$ &&
\\ $Ru$ & $1$ ($K^\circ$) & $2$ (const) & $1$ ($K^\circ$)$^*$ &&&
\\ $HN$ & $1$ ($K^\circ$) & $1$ ($K^\circ$) & $1$ ($K^\circ$) &&&
\\ $ON$ & $1$ (SN) & $1$ ($K^\circ$)$^*$ && $1$ ($K^\circ$) &&
\\ $Ly$ & $1$ (SN) & $1$ ($K^\circ$) & $1$ ($K^\circ$) &&&
\\ $Co_3$ & $1$ (SN) & $1$ ($K^\circ$) & $4$ (nnc)$^*$ &&&
\\ $Co_2$ & $1$ (SN) & $1$ ($K^\circ$) & $S_3$ (nnc)$^*$ &&&
\\ $Co_1$ & $1$ (SN) & $1$ ($K^\circ$) & $1$ ($K^\circ$)$^*$ & $1$ ($K^\circ$)$^*$ &&
\\ $Fi_{22}$ & $1$ (SN) & $1$ ($K^\circ$) & $S_3$ (r2) &&&
\\ $i\cdot Fi_{22}$ & $-$ & $-$ & $i\times S_3$ &&&
\\ $Fi_{23}$ & $1$ (SN) & $1$ ($K^\circ$)$^*$ & $2$ (r2) &&&
\\ $Fi_{24}'$ & $1$ (SN) & $1$ ($K^\circ$)$^*$ & $1$ ($K^\circ$) & $1$ ($K^\circ$)$^*$ && 
\\ $Th$ & $1$ (SN) & $1$ ($K^\circ$)$^*$ & $1$ ($K^\circ$)$^*$ & $2$ (r2)$^*$ &&
\\ $BM$ & $1$ (SN) & $1$ ($K^\circ$)$^*$ & $1$ ($K^\circ$)$^*$ & $2$ (r2)$^*$ &&
\\ $M$ & $1$ (SN) & $1$ ($K^\circ$) & $1$ ($K^\circ$) & $1$ ($K^\circ$) & $1$ ($K^\circ$) & $1$ ($K^\circ$)
\\ $T$ & $1$ (SN) & $1$ ($K^\circ$) & $4A_4$ (TI)&&&
\\ \hline
\end{tabular}
\end{center}
\caption{The groups $N_G(S)/K_G$ (see Section \ref{sec:prelim}) for $G$ a sporadic simple group and $S$ non-cyclic. For $T(G,S)$, take the abelianization of this group. The labels (SN) and so on are indicators of the proof, and are explained in Section \ref{sec:detofgroups}. For $\hat G$ a central extension of $G$, we have $N_{\hat G}(\hat S)/K_{\hat G}\cong N_G(S)/K_G$ unless this table states otherwise. For these lines of the table, a $-$ indicates that $N_{\hat G}(\hat S)/K_{\hat G}\cong N_G(S)/K_G$. $^*$ indicates that the result does not appear elsewhere in the literature.}
\label{tab:result}\end{table}
(The result for $3\cdot J_3$ in characteristic $2$ differs from that in \cite{lassueurmazza2015}. There is a torsion endotrivial module for $M_{11}$ and $p=2$, but it is not trivial source (so does not occur in Table \ref{tab:result}), and comes from the torsion element in $T(S)$.) We must also determine the actions of outer automorphisms on the groups in Table \ref{tab:result}, to understand the group $T(G,S)$ when $G$ is an almost simple sporadic group. This information is given in Section \ref{sec:autos}.

We also attempt to construct all modules in $T(G,S)$ for $G$ a quasisimple sporadic group as well. For $G$ `small' we can do so, but for `large' groups, namely $McL$, $Co_3$, $Co_2$, $J_4$, $Th$ and $BM$, we cannot explicitly construct all of the (non-trivial) elements of $T(G,S)$.

We also consider the case of simple endotrivial modules for quasisimple sporadic groups $G$. In \cite[Table 6]{lassueurmalleschulte2016}, Lassueur--Malle--Schulte give a partial classification of simple endotrivial modules for quasisimple sporadic groups, with (up to duality) fourteen potential simple endotrivial modules undecided. Here we prove or disprove endotriviality for eleven of the fourteen modules left over in \cite{lassueurmalleschulte2016}, four of which were already determined in \cite{lassueurmalle2015} (although we provide alternate proofs for those four modules here).

The proof of Theorem \ref{thm:mainthm} uses a local-group-theoretic characterization of $T(G,S)$ found by Grodal in \cite{grodal2016un}. Because that paper is long and is mainly homotopy theory, in Section \ref{sec:equiv} we give a short proof that shows that Grodal's criterion is equivalent to Balmer's characterization of $T(G,S)$ as the group of `weak homomorphisms' \cite{balmer2013}. The result \cite{balmer2013} uses representation theory, and so together this removes the homotopy theory from the proof.

(The proof given here can be pieced together from various places in \cite{grodal2016un} as well, and of course the homotopy theory is how the characterization was originally found. But it does mean that algebraists have a complete proof of all results in the field.)

Grodal's characterization, given in the next section, turns the computation of $T(G,S)$ into a question about normalizers, which means finite group theory can be used to attack it. This enables us to obtain much better results much more easily than using previous methods; some of the work from several papers \cite{carlsonhemmermazza2010,carlsonmazzanakano2009,lassueurmazza2015a,lassueurmazza2015} can been condensed into just a few pages.

\bigskip

\noindent \textbf{Acknowledgements} I would like to thank: Jesper Grodal for some discussions about his at-the-time unpublished work, and for alerting me to some examples that he and his collaborators have found; Thomas Breuer and Richard Parker for help with computations with large modules for sporadic groups, whose contributions are detailed in the proof of Proposition \ref{prop:knownsimple}; Gunter Malle for reading an early version of this manuscript; the referee for several helpful comments that have improved the paper's exposition.

\section{Preliminary results}
\label{sec:prelim}

Throughout this work $G$ is a finite group, $p$ is a prime and $k$ is an algebraically closed field of characteristic $p$. Let $S$ denote a Sylow $p$-subgroup of $G$, and write $N$ as shorthand for $N_G(S)$. Write $G'$ for the derived subgroup of $G$, and $G^{\mathrm{ab}}$ for $G/G'$, the abelianization of $G$. Our notation for simple groups follows \cite{atlas}, and the notation for conjugacy classes as well. The notation is a number followed by a letter, e.g., 2B, where the number denotes the order of an element and the letter denotes the position of the class among all classes of elements of order $2$, when ordered according to increasing class size. The structure of $N_G(S)$ for $G$ sporadic is given in \cite{wilson1998}.

If $H$ is a subgroup of $G$ write $(-){\downarrow_H}$ for restriction to $H$. Write $V^*$ for the dual of a $kG$-module $V$, and $V^\pm$ to mean $V$ and/or $V^*$. When writing socle structures of a module, we use `/' to delineate the layers, so that
\[ A/B,C/D,E\] 
is a module with socle $D\oplus E$, second socle $B\oplus C$ and third socle $A$.

A $kG$-module $V$ is \emph{endotrivial} if $V\otimes V^*$ is the sum of the trivial module $k$ and a projective module. Every endotrivial module is the sum of a projective module and an indecomposable module. The collection of endotrivial modules forms a group under tensor product if one ignores projective factors, and this group is  denoted $T(G)$. Write $T(G,S)$ for the subgroup of $T(G)$ of all trivial-source modules. Thus $V\in T(G)$ lies in $T(G,S)$ if and only if $V{\downarrow_S}$ is the sum of a single trivial module and a projective module. Note that if $T(S)$ is torsion-free, and this is always the case unless $S$ is cyclic, quaternion or semidihedral \cite{carlsonthevenaz2004}, then $V$ lies in $T(G,S)$ whenever $V$ is self-dual.

There is already a large literature on endotrivial modules for finite groups. We just need a few results from this, mainly for confirming that a given module actually is endotrivial.

\begin{lem}\label{lem:restrict} Let $Q_1,\dots,Q_n$ be representatives from the conjugacy classes of maximal elementary abelian $p$-subgroups of $G$ (i.e., maximal amongst all elementary abelian subgroups). Suppose that the $p$-rank of $G$ is at least $2$. A $kG$-module $V$ is endotrivial if and only if $V{\downarrow_{Q_i}}$ is endotrivial for all $i$. Furthermore, if $G$ is a $p$-group of $p$-rank $2$ then there is a homomorphism
\[ T(G)\to \prod_i T(Q_i)\cong \Z^n.\]
Both the kernel and cokernel of this map are finite. If the $p$-group $G$ is not a semidihedral $2$-group then the kernel is trivial.
\end{lem}

The first statement is \cite[Proposition 2.3(b)]{carlsonmazzanakano2006}. That the kernel of the map is finite is \cite[Proposition 2.3(c)]{carlsonmazzanakano2006}, and that the cokernel is finite follows from the fact that the kernel is finite and the rank of $T(G)$ is equal to $n$ by \cite[Theorem 4]{alperin2001a}. The last statement follows from \cite[Theorem 1.2]{carlsonthevenaz2005}.

We now define a particular subgroup for any finite group $G$ and any prime $p$.

\begin{defn}\label{defn:KandKcirc} Let $G$ be a finite group with Sylow $p$-subgroup $S$. The subgroup $K_G$ of $G$ is defined as consisting of all elements $g\in N_G(S)$ such that there exist $x_1,\dots,x_n,y\in G$ and subgroups $Q_1,\dots,Q_n$ of $S$, with the following properties:
\begin{enumerate}
\item $g=x_1\dots x_n$;
\item $x_i\in O^{p'}(N_G(Q_i))$ for all $1\leq i\leq n$;
\item $y^{x_1\dots x_i}$ lies in both $Q_i$ and $Q_{i+1}$, for all $1\leq i<n$.
\end{enumerate}
Let $K_G^\circ$ be the normal subgroup of $N_G(S)$ generated by $N_G(S)\cap O^{p'}(N_G(Q))$ for all $N_G(S)$-conjugacy classes of subgroups $1<Q\leq S$.
\end{defn}

It is clear that $K_G^\circ\leq K_G\normal N_G(S)$. Often $K_G=K_G^\circ$, but not always. The main result of Grodal's is the following.

\begin{thm}[{{Grodal \cite[Theorem 4.27]{grodal2016un}}}]\label{thm:grodal} There is an isomorphism
\[ (N_G(S)/K_G)^{\mathrm{ab}}\cong T(G,S).\]
\end{thm}

We therefore must compute $K_G$ for $G$ a (quasi)simple group. It is clear that for a specific group, if one can compute $O^{p'}(N_G(Q))$, say on a computer, then one can easily compute $K_G^\circ$. The computation of $K_G$ looks more difficult, but luckily in all but one case that we consider here $K_G^\circ=K_G$. (In fact, the author knows of only one other case where $K_G\neq K_G^\circ$, that of $G_2(5)$ and $p=3$, to which he was alerted by Grodal. It seems very rare for quasisimple groups that $K_G$ differs from $K_G^\circ$.)

If a Sylow $p$-subgroup $S$ of $G$ is cyclic then we see that $K_G=S$ since $O^{p'}(N_G(S))=S$, so we obtain that $T(G,S)$ is isomorphic to the abelianization of $N_G(S)/S$, agreeing with \cite{mazzathevenaz2007}. We therefore concern ourselves with groups with non-cyclic Sylow $p$-subgroups.

Another consequence of Grodal's result is the following previously known statement \cite[Proposition 2.8]{carlsonmazzanakano2006}.
\begin{lem}\label{lem:ti} If $G$ has trivial intersection Sylow $p$-subgroups then $K_G=S$, so $T(G,S)\cong (N_G(S)/S)^{\mathrm{ab}}$.
\end{lem}

%
%
%
%
%

If $K_G^\circ=N_G(S)$ then obviously $K_G=N_G(S)$ and we are done. However, if $K_G^\circ\neq N_G(S)$, then we need a criterion to guarantee $K_G^\circ=K_G$; the next lemma furnishes us with such a criterion.

\begin{lem}\label{lem:nonon-cyclic} If $O^{p'}(N_G(Q))\leq S$ for all $1<Q<S$ not of order $p$, then $K_G=K_G^\circ$. In particular, if $|S|=p^2$ then $K_G=K_G^\circ$.
\end{lem}
\begin{pf} From the assumption, the subgroups $Q_i$ in the definition of $K_G$ must all have order $p$, or else the $x_i$ all lie in $S$. However, if $Q_i$ has order $p$ then $Q_i=\gen y$ for all $i$, and so $g\in O^{p'}(N_G(Q))\leq K_G^\circ$. Thus $K_G\leq K_G^\circ$, as claimed.
\end{pf}

We illustrate this with an example, solving some of the open cases from \cite{lassueurmazza2015}.

\begin{example} Let $G=J_3$ and $p=3$. The order of $S$ is $3^5$, but the only subgroup $Q$ of $S$ for which $O^{p'}(N_G(Q))$ is not contained in $S$ is of order $3$. Thus $K_G=K_G^\circ$ by Lemma \ref{lem:nonon-cyclic}, and in this case $N_G(S)/K_G^\circ$ has order $2$ by a direct computation on Magma \cite{magma}. We may also apply Lemma \ref{lem:nonon-cyclic} to $Co_2$ and $Co_3$ for $p=5$ (in both cases the size of $S$ is $5^3$), and in these situations $N_G(S)/K_G$ is isomorphic to $C_4$ and $S_3$ respectively.
\end{example}

We will also want a simple criterion to guarantee that a given element of $N_G(S)$ lies in $K_G^\circ$. Note that, trivially, if $Q$ has order $p$, $O^{p'}(N_G(Q))\leq C_G(Q)$. Often centralizers of $p$-elements in big sporadic groups satisfy $O^{p'}(C_G(z))=C_G(z)$, so we can just test if $g\in N_G(S)$ lies in $C_G(z)$ for some $z\in S$ with this property.

\begin{lem}\label{lem:autocent} Let $z$ be an element of order $p$ in $G$ such that $O^{p'}(C_G(z))=C_G(z)$, and let $g\in N_G(S)$. If $C_S(g)$ contains a $G$-conjugate of $z$ then $g\in K_G^\circ$.
\end{lem}

We give an example of how to use this lemma now, dealing with another open case from \cite{lassueurmazza2015}.

\begin{example}\label{ex:thchar3} Let $G=Th$, $p=3$. The centralizer of a 3A element is $3\times G_2(3)$, so satisfies the condition of Lemma \ref{lem:autocent}. The group $N_G(S)/S$ is isomorphic to the Klein four group. Let $E$ be a Klein four complement to $S$ in $N_G(S)$.

To prove that every element of $E$ centralizes a 3A element, we check on a computer. We start with a $248$-dimensional representation of $G$ from the online Atlas \cite{webatlas}, and a straight-line program that gives us one of the two $3$-local subgroups that contain $N_G(S)$. Let \texttt{H} be this group, with Sylow $3$-subgroup \texttt{S}. To more easily compute inside \texttt{H}, it is better to produce a copy \texttt{H2} of \texttt{H} as a power-commutator (or equivalently power-conjugate) group. We construct the subgroup \texttt{E} inside \texttt{H2}, then compute the centralizer in \texttt{S} of each non-trivial element of \texttt{E} (mapped back from \texttt{H2} into \texttt{H}) and check that they contain a 3A element.

Inside \texttt{H} we can detect the conjugacy classes of elements of order $3$ by the size of their fixed-point space on the representation of $G$. We provide Magma code for the whole procedure, starting from the point where we assume we have \texttt{H}:
\begin{verbatim}
> H2,phi:=PCGroup(H); S:=Sylow(H2,3); E:=Sylow(Normalizer(H2,S),2);
> Order(E);
4
> for x in [E.1,E.2,E.1*E.2] do
for> C:=Centralizer(S,x) @@phi;
for> {*Dimension(Eigenspace(i,1)):i in C|Order(i) eq 3*};
for> end for;
{* 80^^6, 86^^26, 92^^12 *}
{* 80^^6, 86^^26, 92^^12 *}
{* 80^^6, 86^^26, 92^^12 *}
\end{verbatim}
This shows that all three of the classes of elements of order $3$ lie in the fixed-point space of each involution in $E$. In particular, $N_G(S)=K_G^\circ$.

In fact, we see from \cite{wilson1988} that $O^{3'}(C_G(z))=C_G(z)$ for all elements $z$ of order $3$, so we would just need to check that every element of $E$ has a fixed point on $S$, but for most groups not all centralizers of elements of order $p$ have this property, so this extra check needs to be implemented.
\end{example}

With these two tricks, the information on local subgroups from the Atlas \cite{atlas} and some of their constructions on the online Atlas \cite{webatlas}, and the structures of $N_G(S)/S$ given in \cite{wilson1998}, we can fairly quickly despatch all groups in Section \ref{sec:detofgroups}.

\section{Equivalence of Grodal's and Balmer's characterizations}
\label{sec:equiv}

In \cite{balmer2013}, Balmer gives a group-theoretic criterion for the determination of $T(G,S)$, with an algebraic proof. He defines a group of `weak $H$-homomorphisms' $A(G,H)$, and for $H=S$ proves that $A(G,S)\cong T(G,S)$. Here we assume throughout that $H=S$, as we have no need for the more general definition. However, we want to generalize weak $S$-homomorphisms in a different way.

\begin{defn}
A function $\rho:G\to \GL_n(k)$ is a \emph{weak $n$-homomorphism} if $\rho(x)=1$ (i.e., the identity matrix) whenever $S\cap S^x=1$ or $x$ is a $p$-element normalizing some non-trivial subgroup of $S$, and if $S\cap S^y\cap S^{xy}\neq 1$ then $\rho(x)\rho(y)=\rho(xy)$.
\end{defn}

Let $A_n(G,S)$ denote the set of weak $n$-homomorphisms $G\to \GL_n(k)$, and $A(G,S)=A_1(G,S)$. We will prove the following.

\begin{thm}\label{thm:mainresult} The restriction map from $G$ to $N_G(S)$ induces a bijection
\[ A_n(G,S)\to \Hom(N_G(S)/K_G,\GL_n(k)).\]
\end{thm}

The case $n=1$ shows that Balmer's and Grodal's characterizations of $T(G,S)$ are the same. For $n\geq 1$, Grodal showed \cite[Theorem 3.10]{grodal2016un} that every indecomposable element of $\Hom(N_G(S)/K_G,\GL_n(k))$ is the Green correspondent of an indecomposable $kG$-module, whose restriction to $S$ is the sum of a free and a trivial module (of some dimension). For $n\geq 2$ the content of Theorem \ref{thm:mainresult} is to give the correct extension of the definition of weak $S$-homomorphisms from \cite{balmer2013} to arbitrary dimensions.

\medskip
We first notice that if $\rho$ is a weak $n$-homomorphism and $1<Q\leq S$, then the restriction of $\rho$ to $N_G(Q)$ is a homomorphism to $\GL_n(k)$. In particular, the restriction of any weak $n$-homomorphism to $N_G(S)$ is a map in $\Hom(N_G(S),\GL_n(k))$. Moreover, since $x\in N_G(Q)$ is a $p$-element then $\rho(x)=1$, so therefore $O^{p'}(N_G(Q))$ lies in the kernel of this homomorphism.

(Comparing our definition to the original definition of weak $1$-homomorphism from \cite{balmer2013}, this extra condition on $p$-elements is not present. However, since $k^\times$ has no $p$th roots of unity, if $x\in N_G(Q)$ is a $p$-element then $\rho(x)=1$ for any $\rho\in A_1(G,S)$---since $\rho(x^i)=\rho(x)^i$ by the weak $n$-homomorphism property---so the condition always holds for weak $S$-homomorphisms in the sense of \cite{balmer2013}.)

Let $g\in G$ and $A$ be a non-trivial subset of $S\cap S^{g^{-1}}$. An \emph{Alperin decomposition} of $g$ with respect to $A$ (from Alperin's fusion theorem \cite{alperin1967}) is an expression $g=x_1\ldots x_rz$, where $x_i\in N_G(Q_i)$ for some subgroup $1<Q_i\leq S$, $z\in N_G(S)$, each $x_i$ is a $p$-element, and $A^{x_1\ldots x_{i-1}}\leq Q_i$ for all $1\leq i<r$.

\begin{lem}\label{lem:basicprops} Let $\rho:G\to \GL_n(k)$ be a weak $n$-homomorphism.
\begin{enumerate}
\item If $g=x_1\ldots x_rz$ is an Alperin decomposition of $g$ with respect to any subset non-trivial $A$, then $\rho(g)=\rho(z)$.
\item  If $\rho(g)=1$ for all $g\in N_G(S)$ then $\rho(g)=1$ for all $g\in G$.
\item The restriction map from $A_n(G,S)$ to $\Hom(N_G(S),\GL_n(k))$ is injective.
\end{enumerate}
\end{lem}
\begin{pf} Suppose we have an Alperin decomposition $g=x_1\ldots x_rz$. Since $x_i$ is a $p$-element of $N_G(Q_i)$ for $1<Q_i\leq S$, we have that $\rho(x)=1$. We claim that the weak $n$-homomorphism property implies that $\rho(x_1\ldots x_i)=1$ for all $1\leq i\leq r$, and then that $\rho(g)=\rho(z)$.

We proceed by induction on $i$, the result obviously holding that $\rho(x_1)=1$. Thus write $x=x_1\ldots x_{i-1}$, $x'=x_i$ and suppose that $\rho(x)=1$ (and $\rho(x')=1$). We merely need to show that
\[ S\cap S^{x'}\cap S^{xx'}\neq 1,\]
but this is equivalent to
\[ S^{x'^{-1}}\cap S\cap S^{x}\neq 1;\]
the subgroup $R^x$ lies in $S^{x'^{-1}}$, $S$, and $S^x$, so this triple intersection is indeed non-empty. Thus $\rho(xx')=\rho(x)\rho(x')$, and thus $\rho(x_1\ldots x_r)=1$. Finally, since $z\in N_G(S)$, we clearly have that the triple intersection property holds for $x_1\ldots x_r$ and $z$, and therefore $\rho(g)=1\cdot \rho(z)$, proving (i).

\medskip

We now prove (ii), so let $g\in G$. If $S\cap S^g=1$ then $\rho(g)=1$ by definition. If $S\cap S^g\neq 1$ then we obtain an Alperin decomposition of $g$ with respect to $S\cap S^{g^{-1}}$, and so $\rho(g)=\rho(z)$ for some $z\in N_G(S)$. However, $\rho(z)=1$ by assumption, so $\rho(g)=1$, as claimed.

\medskip

We now prove (iii). Since the product of two weak $n$-homomorphisms is a weak $n$-homomorphism, if $\rho\sigma^{-1}$ is trivial on $N_G(S)$ then $\rho\sigma^{-1}=n\cdot 1_G$, but then $\rho=\sigma$, as claimed.
\end{pf}

\begin{lem} Let $\rho:G\to \GL_n(k)$ be a weak $n$-homomorphism. If $g\in K_G$ then $\rho(g)=1$. In particular, the restriction map yields an injection $A_n(G,S)\to \Hom(N_G(S)/K_G,\GL_n(k))$.
\end{lem}
\begin{pf} By definition, $K_G$ consists of those elements $g$ of $N_G(S)$ with an Alperin decomposition $g=x_1\ldots x_rz$ with $z=1$. By Lemma \ref{lem:basicprops}(i), $\rho(g)=\rho(z)=1$, and so $K_G\leq \ker (\rho{\downarrow_{N_G(S)}})$. This completes the proof.
\end{pf}

We now must construct a map going the other direction. Thus, given any homomorphism $\rho:N_G(S)\to \GL_n(k)$ with $K_G$ in the kernel, we need to extend this to a weak $n$-homomorphism $G\to \GL_n(k)$ in $A_n(G,S)$. We do this as follows.
\begin{enumerate}
\item If $g\in G$ satisfies $S\cap S^g=1$, then set $\rho(g)=1$.
\item Thus $S\cap S^g\neq 1$, so we have an Alperin decomposition $g=x_1\ldots x_rz$. Set $\rho(g)=\rho(z)$.
\end{enumerate}

Call this the \emph{weak $n$-extension} of the map $\rho:N_G(S)\to \GL_n(k)$. We claim it is well defined: let $x_1\ldots x_rz$ be an Alperin decomposition of $g$ with respect to $S\cap S^{g^{-1}}$, and let $x_1'\ldots x_s'z'$ be any other Alperin decomposition of $g$, with respect to some (non-trivial) subset $A'$ of $S$. We have that
\[ x_1\ldots x_rz=x_1'\ldots x_s'z',\]
so
\[ zz'^{-1}=x_r^{-1}\ldots x_1^{-1}x_1'\ldots x_s'.\]
We claim that the right-hand side is an Alperin decomposition of $zz'^{-1}$, and therefore that $zz'^{-1}$ lies in $K_G$. Thus $\rho(z)=\rho(z')$, as needed.

We now check that it is indeed an Alperin decomposition. To see this, we set $A=A'^g\leq S$. It is easy to see that this expression satisfies the condition of being an Alperin decomposition with respect to $A$, so we see that $\rho$ is well defined.

\begin{prop} If $\rho:N_G(S)\to \GL_n(k)$ is a homomorphism with $K_G$ in the kernel, then the weak $n$-extension (also denoted $\rho$) to $G$ is a weak $n$-homomorphism. Consequently, the restriction map induces an isomorphism
\[ A_n(G,S)\to \Hom(N_G(S)/K_G,\GL_n(k)),\]
and so in particular $T(G,S)$ consists of the Green correspondents in $G$ of elements of $\Hom(N_G(S)/K_G,k^\times)$.
\end{prop}
\begin{pf} Let $g\in G$. If $g\in S$ or $S\cap S^g=1$ then $\rho(g)=1$ by construction. Thus we need to check that $\rho(gh)=\rho(g)\rho(h)$ whenever $S\cap S^h\cap S^{gh}\neq 1$. Note that, in particular each of $S\cap S^h$, $S\cap S^{gh}$, and $S\cap S^g$ is non-trivial. Let
\[ g=x_1\ldots x_rz,\quad h=x_1'\ldots x_s'z',\]
be Alperin decompositions of $g$ and $h$, with respect to the subsets $S\cap S^{g^{-1}}$ and $S\cap S^{h^{-1}}$ respectively. We claim that
\begin{equation} gh=\left[x_1\ldots x_r(x_1'\ldots x_s')^{z^{-1}}\right]\left[zz'\right]\label{eq:alpdec}\end{equation}
is an Alperin decomposition of $gh$, and therefore 
\[\rho(gh)=\rho(zz')=\rho(z)\rho(z')=\rho(g)\rho(h),\]
since $z,z'\in N_G(S)$.

Certainly $gh$ is equal to the expression, so we need to find a subset $A$ that is transported by each $x_i$ and $(x_j')^{z^{-1}}$. Set $A=S^{h^{-1}g^{-1}}\cap S^{g^{-1}}\cap S\neq 1$. Let $x_i\in N_G(Q_i)$ and $x_i'\in N_G(Q_i')$ for each $i$: we have that $A\leq S\cap S^{g^{-1}}$, so that $A^{x_1\ldots x_{i-1}}\leq Q_i$ by the definition of an Alperin decomposition.

Set $B=A^{x_1\ldots x_r}=A^{gz^{-1}}$. Notice that since $A^g\leq Q_1'$, $B\in (Q_1')^{z^{-1}}$, and indeed we see that $B^{x_1'\ldots x_{i-1}'}\leq (Q_i')^{z^{-1}}$ for all $1\leq i<s$. Thus the expression in (\ref{eq:alpdec}) is an Alperin decomposition, as needed.
\end{pf}

\section{Determination of $K_G^\circ$ and $K_G$}
\label{sec:detofgroups}

If $G$ is, say, $M_{11}$, one may load the group directly into Magma and in a few seconds compute $K_G^\circ$. For groups where there is a permutation representation available, and Magma can determine $K_G^\circ$ in less than a couple of minutes, we omit a proof that $K_G^\circ$ is as claimed. In all cases except for $G=Ru$ for $p=3$ and $G=3\cdot J_3$ for $p=2$, $K_G=K_G^\circ$ either by an application of Lemma \ref{lem:nonon-cyclic} or because $K_G^\circ=N_G(S)$. For $Ru$ we construct an endotrivial module to demonstrate that $T(G,S)$ has order at least $2$, and for $3\cdot J_3$ we express an element of $K_G\setminus K_G^\circ$ as the product of two elements $x_1x_2$, as in Definition \ref{defn:KandKcirc}.

The explanations in Table \ref{tab:result} in brackets are shorthand, to alert the reader as to which technique is used to conclude that $K_G$ is as stated, rather than just $K_G^\circ$. (SN) means that $G$ has a self-normalizing Sylow $2$-subgroup, so $N_G(S)=S$. ($K^\circ$) means $K_G^\circ$ is already all of $N_G(S)$. (r2) means $S\cong C_p\times C_p$, so we may apply Lemma \ref{lem:nonon-cyclic} above to see that $K_G=K_G^\circ$. (nnc) means that for all $Q\leq S$ with $Q$ not cyclic of order $p$, $O^{p'}(N_G(Q))\leq S$, hence $K_G=K_G^\circ$ by Lemma \ref{lem:nonon-cyclic} above. (TI) means that $G$ has a trivial intersection Sylow $p$-subgroup, so that $K_G=S$ by Lemma \ref{lem:ti}. Finally, (const) means that $|N_G(S):K_G^\circ|$ was found to be equal to $2$, and a trivial-source endotrivial module other than $k$ was constructed directly.

This section proves the following result as a consequence of a case-by-case analysis.

\begin{thm}\label{thm:sporadic} If $G$ is a sporadic quasisimple group and $S$ is non-cyclic then the group $N_G(S)/K_G$ is given in Table \ref{tab:result}. In addition, we have $K_G=K_G^\circ$, unless $G=3\cdot J_3$ and $p=2$, in which case $|K_G:K_G^\circ|=3$. Moreover, in this case $N_G(S)=K_G$.
\end{thm}

We first compute $K_G$ in the case where $G$ is simple in Section \ref{sec:simple}. In Section \ref{sec:central} we determine for $G$ quasisimple but not simple, whether $K_G$ contains $Z(G)$. These two subsections furnish us with a proof of Theorem \ref{thm:sporadic}.

\subsection{Simple groups}
\label{sec:simple}
We begin by proving Theorem \ref{thm:sporadic} for the sporadic simple groups $G$. We will give example Magma code in a few of the earlier cases to show how to compute the information we use.

\subsubsection{The groups $M_{11}$, $M_{12}$, $M_{22}$, $M_{23}$, $M_{24}$, $J_1$, $J_2$, $HS$, $McL$, $He$, $Co_3$, $Co_2$, $T$}

In all cases there is a small permutation representation of $G$ available via \cite{webatlas} and one may easily determine $K_G^\circ$; enumerate all conjugacy classes of subgroups $Q$ of $S$, take all normalizers $N_G(Q)$ of them, and then take the subgroup generated by all normal closures of $S$ inside all $N_G(Q)$, to produce $K_G^\circ$. It is as in Table \ref{tab:result}.

\subsubsection{The group $Ru$ $(p=2,3,5)$}

Unless $p=3$, an easy computer calculation determines that $K_G^\circ=N_G(S)$. For $p=3$, we have that $|N_G(S):K_G^\circ|=2$, and we cannot apply Lemma \ref{lem:nonon-cyclic}. However, there is a simple endotrivial module with trivial source, namely the $406$-dimensional simple module \cite{lassueurmalleschulte2016}. One may explicitly check that it has trivial source, or note that it is self-dual and $T(S)$ is torsion-free. Thus $T(G,S)\cong C_2$, and in particular $K_G=K_G^\circ$.

\subsubsection{The group $Suz$ $(p=2,3,5)$}

If $p$ is odd then an easy computer calculation determines $K_G^\circ$. If $p=2$ then $N_G(S)\cong S\rtimes 3$ (see \cite{wilson1998}). There are many $2$-subgroups, so we give a quick theoretical proof that $K_G^\circ=N_G(S)$. If $z$ is an involution in $G$ then $C_G(z)$ has no non-trivial odd quotients as the two centralizers have composition factors only simple groups and $C_2$ (see \cite[Table 5.3o, p.276]{glsvol3}). Thus it suffices to check that if $E$ denotes a Sylow $3$-subgroup of $N_G(S)$, then $C_S(E)\neq 1$. This is the case since $|S|=2^{13}$ and so $3\nmid (|S|-1)$. Thus $K_G^\circ=N_G(S)$.

\subsubsection{The group $HN$ $(p=2,3,5)$}

Here there is a permutation representation, but it is large, so we give more details to prove Theorem \ref{thm:sporadic} quickly, again assisted by a computer to reduce errors. For $p=2$, the same method of proof as for $Suz$ applies. This time the table is \cite[Table 5.3w, p.284]{glsvol3}, but as $|S|=2^{14}$, in theory we could have $C_S(E)=1$. In fact, $|C_S(E)|=64$ by a Magma calculation:
\begin{verbatim}
> S:=Sylow(G,2); //G is a copy of HN
> NGS:=Normalizer(G,S); E:=Sylow(NGS,3);
> Order(Centralizer(S,E));
64
\end{verbatim}

Let $p=3$. Here $|S|=3^6$, $|N_G(S)/S|=8$. All elements $z$ of order $3$ in $G$ satisfy $O^{3'}(C_G(z))=C_G(z)$ (again, see \cite[Table 5.3w, p.284]{glsvol3}), so $K_G^\circ$ contains $N_G(S)\cap C_G(z)=C_N(z)$ (where $N=N_G(S)$) for all elements $z\in \Omega_1(S)\setminus \{1\}$. A quick computer calculation shows that every element of a Sylow $2$-subgroup of $N$ lies in $C_N(z)$ for some $z$ of order $3$, and therefore $K_G^\circ=N$.
\begin{verbatim}
> S:=Sylow(G,3); //G is a copy of HN
> NGS:=Normalizer(G,S); E:=Sylow(NGS,2);
> [Order(Centralizer(S,i)):i in E];
[ 729, 3, 9, 3, 9, 3, 9, 3 ]
\end{verbatim}

Now let $p=5$. Here $|S|=5^6$, $|N_G(S)/S|=8$. There are three classes of elements (classes 5B, 5C and 5D) that are fixed-point-free in the permutation representation on 1140000 points. For $z$ from each of these classes, $O^{5'}(C_G(z))=C_G(z)$, so $K_G^\circ$ contains the subgroup of $N=N_G(S)$ generated by $C_N(z)$ for all fixed-point-free $z\in S$. This is all of $N_G(S)$, as in the $p=3$ case.

\subsubsection{The group $Co_1$ $(p=2,3,5,7)$}
If $p=2$ then $N_G(S)=S$. If $p=5,7$, then $K_G^\circ=N_G(S)$ and the computation is easy on a computer. For $p=3$ it is much more difficult because the Sylow $3$-subgroup is large enough to make enumerating all subgroups challenging.

In this case, $|S|=3^9$ and $|N_G(S)/S|=8$. Let $z_1$ be an element of $S$ with centralizer $3\cdot Suz$, and $z_2$ be an element of $S$ with centralizer $3\times A_9$. The $G$-conjugates of $z_1$ in $S$ split into two $N_G(S)$-orbits, and the $G$-conjugates of $z_2$ in $S$ split into five $N_G(S)$-orbits. Then $N=N_G(S)$ is generated by $C_N(z_1)$ and $C_N(z_2)$, unless $z_2$ lies in the unique class with $N_G(S)$-centralizer of order $243$. Thus if $z_2$ is chosen not to lie in this case, $N$ is generated by $C_N(z_1)$ and $C_N(z_2)$, both of which lie in $K_G^\circ$. Hence $K_G^\circ=N_G(S)$, as claimed.

\subsubsection{The groups $Fi_{22}$, $Fi_{23}$ and $Fi_{24}'$ $(p=2,3,5,7)$}
In each of these cases a permutation representation is available, and can easily be used to determine $K_G^\circ$ for $p=5,7$. (For $p=7$ we have $G=Fi_{24}'$ and $K_G^\circ=N_G(S)$, and if $p=5$ then $S\cong C_p\times C_p$, so we may apply Lemma \ref{lem:nonon-cyclic}.) If $p=2$ then $N_G(S)=S$ so we are done, and therefore only $p=3$ remains. As with $Co_1$, there are too many $3$-subgroups to easily enumerate all of them, so we give more details.

If $G=Fi_{22}$, then $|S|=3^9$ and $|N_G(S)/S|=4$. Let $z$ be an element from class 3B, whose centralizer has no $3'$-quotients. There are thirteen $N_G(S)$-classes of 3B elements $x_i$, and together the $C_N(x_i)$ generate $N_G(S)$.

If $G=Fi_{23}$, then $|S|=3^{13}$ and $|N_G(S)/S|=8$. Let $z$ be an element from class 3A, whose centralizer is $3\times O_7(3)$. There are seven $N_G(S)$-classes with 3A elements $x_i$, and together the $C_N(x_i)$ generate $N_G(S)$.

If $G=Fi_{24}'$, then $|S|=3^{16}$ and $|N_G(S)/S|=8$. All elements $z$ of order $3$ in $G$ satisfy $O^{3'}(C_G(z))=C_G(z)$ (see \cite[p.82]{wilson1987}), so $K_G^\circ$ contains the subgroup of $N=N_G(S)$ generated by $C_N(z)$ for all $z\in S$. This is all of $N_G(S)$. One does not in fact need all elements. An element from class 3A has centralizer $3\times O_7(3)$. There are seven $N_G(S)$-classes with 3A elements $x_i$, and together the $C_N(x_i)$ generate $N_G(S)$.

\medskip

This is the end of the list of groups where there is a permutation representation available, and hence one can be very explicit. The rest only have matrix representations available, and so can require slightly more theoretical proofs, for the most part.

\subsubsection{The group $J_4$ $(p=2,3,11)$}
If $p=2$ then $N_G(S)=S$, and if $p=11$ then $S$ is TI, so we may apply Lemma \ref{lem:ti}. Thus $p=3$.

Here $|S|=3^3$, $|N_G(S)/S|=32$. Let $E$ denote a Sylow $2$-subgroup of $N_G(S)$. The structure of $N_G(S)$ is given in \cite[p.302]{wilson1998}; it is $(2\times 3^{1+2}\rtimes 8)\rtimes 2$, and lies inside the maximal subgroup $2^{1+12}.3\cdot M_{22}\rtimes 2$, in fact inside the normalizer of an element of order $3$, which is $6\cdot M_{22}.2$. The element of order $3$ therefore has centralizer $6\cdot M_{22}$, so we may apply Lemma \ref{lem:autocent}: if $E$ is generated by elements that are not fixed-point-free on $S$, then $K_G^\circ$ contains $E$. We can check this easily inside $6\cdot M_{22}$, and find $24$ such elements, obviously enough to generate a group of order $32$. (Alternatively but similarly, we can also note that we already find in $K_G^\circ$ a subgroup of index $2$ in $E$, coming from $C_G(Z(S))\cong 6\cdot M_{22}$, so we consider the sixteen elements in $E\setminus C_E(Z(S))$, and half of these centralize an element of order $3$ on $S$.) Thus $K_G^\circ=N_G(S)$, as claimed.

\subsubsection{The group $Th$ $(p=2,3,5,7)$}

If $p=2$ then $N_G(S)=S$, and if $p=3$ then this is Example \ref{ex:thchar3}. Thus $p=5,7$.
 
If $p=7$ then $|S|=7^2$ (so we may apply Lemma \ref{lem:nonon-cyclic} to obtain $K_G=K_G^\circ$) and $|N_G(S)/S|=144$. There is a single class of elements $x$ of order $7$ in $G$, with centralizer $C_G(x)\cong \gen x\times \PSL_2(7)$ (see \cite[Table 5.3x, p.285]{glsvol3}). Thus $K_G^\circ$ is the subgroup generated by $C_N(x)$ for all $x\in S\setminus\{1\}$. This subgroup has order $72\cdot |S|$, and so $N_G(S)/K_G$ has order $2$. (To see this, note that all elements of order $7$ are conjugate in $G$, and therefore we obtain a subgroup of order $3$ acting non-trivially each of the eight subgroups of order $7$ and centralizing a complement. Inside $N_G(S)/S\cong 3\times 2\cdot S_4$, we find three classes of subgroups of order $3$, with one, four and eight subgroups in the class. Thus the subgroups of order $3$ must be the single class of eight, and they generate the subgroup $3\times 2\cdot A_4$, of index $2$ in $N_G(S)/S$.)

For $p=5$, we have $|S|=5^3$ and $|N_G(S)/S|=96$. There is a single $N_G(S)$-class of subgroups $Q$ of order $25$, with normalizer (it is a maximal subgroup of $G$) $5^2\rtimes \GL_2(5)$, so the subgroup we are interested in is $5^2\rtimes \SL_2(5)$. The intersection of this with $N=N_G(S)$ is $S\rtimes 4$, with the $4$ acting as a torus element on $Q$. This determines a unique subgroup $5^2\rtimes 4$ of $N_N(Q)$, and so we may work entirely inside $N$. Taking the normal closure in $N$ of the subgroup $S\rtimes 4$ yields all of $N$, so $K_G=K_G^\circ=N_G(S)$. We give some example Magma code to ease the reader's verification of this.
\begin{verbatim}
// Assume you are given a copy NGS of the normalizer of a Sylow 5-subgroup of Th
S:=Sylow(NGS,5); Q:=Subgroups(NGS:OrderEqual:=25)[1]`subgroup;
a:=Q.1; b:=Q.2; Q eq sub<Q|a,b>; // Check that this is true
NNQ:=Normalizer(NGS,Q);
g:=[i:i in NNQ|a^i eq a^2 and b^i eq b^-2]; #g; //Check that this has order 25
g:=g[1]; A:=sub<NGS|S,g>;
Index(NGS,NormalClosure(NGS,A)); //Should equal 1
\end{verbatim}

\subsubsection{The Group $Ly$ $(p=2,3,5)$}

If $p=2$ then $N_G(S)=S$. If $p=5$ then $N_G(S)$ is contained in a subgroup $5^3.\PSL_3(5)$. This is already a $5$-local subgroup, so $K_G^\circ=N_G(S)$.

Thus $p=3$. We ape the proof in Example \ref{ex:thchar3}, with a slight twist to speed things up. The order of $S$ is $3^7$, and $N_G(S)/S\cong 2\times SD_{16}$. Let $E$ denote a Sylow $2$-subgroup of $N_G(S)$. We see that $N_G(S)$ is contained in a maximal subgroup $3^5\rtimes(2\times M_{11})$, so there is a $3$-local subgroup $3^5\rtimes M_{11}$. This proves that the $SD_{16}$ subgroup $E_1$ of $E$ is contained in $K_G^\circ$, so we simply need to find one element $z$ of $E\setminus E_1$ in $K_G^\circ$ and we are done. Note that the centralizer of a 3A element is $3\cdot McL$, so we may apply Lemma \ref{lem:autocent}.

There are two obvious candidates for $z$: the two central involutions in $E\setminus E_1$. In both cases $C_S(z)$ has order $9$ and contains no element from 3A. However, of the fourteen other elements, ten have 3A elements in $C_S(z)$. (One may tell 3A and 3B apart because they have different fixed points (dimensions $39$ and $21$ for 3A and 3B respectively) on the $111$-dimensional representation over $\F_5$. Six elements from $E\setminus E_1$ have both classes in their centralizer.)

\subsubsection{The group $BM$ $(p=2,3,5,7)$}

If $p=2$ then $N_G(S)=S$. If $p=7$, then $|S|=7^2$, and the structure of $N_G(S)$ is
\[ (2^2\times (S\rtimes (3\times 2A_4))).2.\]
It lies in the maximal subgroup $(2^2\times F_4(2)).2$ \cite[Table II]{wilson1998}, which is the centralizer of a 2C involution. (The outer 2 swaps the two 2A involutions in the centre of the derived subgroup.) There is a single conjugacy class of elements of order $7$, with centralizer of structure $7\times 2.\PSL_3(4)\rtimes 2$, so $O^{7'}(C_G(z))\cong 7\times 2\cdot \PSL_3(4)$. This lies in the centralizer of a 2A involution, so the 2A involution in $C_G(S)$ lies in $K_G^\circ$.

Also, for each subgroup $Q$ of order $7$, $K_G^\circ$ contains an element of order $3$ that centralizes $Q$ and normalizes a complement in $S$. Each of these eight subgroups $3$ must be different, so $K_G^\circ$ contains all of $3\times 2A_4$. Thus $K_G^\circ$ contains all of $2^2\times 3\times 2A_4$, so a subgroup of index $2$ in $N_G(S)$. (This can also be seen as $Th\leq BM$.)

We claim this is all of $K_G^\circ$, and hence $K_G$ by Lemma \ref{lem:nonon-cyclic}. To see this note that this subgroup of index $2$ has two conjugacy classes of subgroups of order $7$, which are interchanged by the top $2$ in $N_G(S)$. Thus this top $2$ does not centralize any subgroup of order $7$, and therefore cannot appear in $C_G(Q)$. Thus $K_G=K_G^\circ$ has index $2$ in $N_G(S)$, as claimed.

Now let $p=5$. By \cite[Table I]{wilson1998}, $N_G(S)/S\cong 4\times 4$. Furthermore, $5^3\cdot \PSL_3(5)$ is a $5$-local subgroup \cite{wilson1987}, and clearly has no $5'$-quotients; it also contains $N_G(S)$, so $K_G^\circ=N_G(S)$, as needed.

Finally, let $p=3$. By \cite[Table I]{wilson1998}, $N_G(S)/S\cong 2^3$. Note that $S$, and indeed $N_G(S)$, is contained in the subgroup $H\cong Fi_{23}$, and since $N_G(S)=K_H^\circ$ as we have seen above, $N_G(S)=K_G^\circ$ as well.

\subsubsection{The Group $M$ $(p=2,3,5,7,11,13)$}

If $p=2$ then $N_G(S)=S$. Suppose that $p=13$. By \cite[Table I]{wilson1998}, $N_G(S)/S\cong 3\times 4S_4$, and $S\cong 13^{1+2}_+$. Since $N_G(S)$ is a maximal subgroup of $G$ (it is the normalizer of a 3B element), and we cannot work inside $G$ directly, we will use other $13$-local subgroups and decide which elements of $N_G(S)$ lie in them. The centralizer of a 13A element is $13\times \PSL_3(3)$, and there is also maximal subgroup $13^2\rtimes \SL_2(13).4$, with the $13^2$ being pure, and hence consisting of 13B elements.

There are two $N_G(S)$-classes of subgroups of order $13^2$, one with an involution in $N_G(S)$ inverting all elements of it and the other without. Thus we have identified a subgroup $Q$ of $S$ with normalizer $13^2\rtimes \SL_2(13).4$. Fix an element $z$ of order $2$ inverting $Q$ (all such elements lie in $Q\gen z$), and we aim to find a Borel subgroup of the $\SL_2(13)$. There are $104$ elements of order $12$ powering to $z$ in $N_S(Q)$, yielding two $N_G(S)$-conjugacy classes of subgroups $(13^2)\rtimes (13\rtimes 12)$. However, the normal closure in $N_G(S)$ of any of the $104$ such subgroups is all of $N_G(S)$, so $K_G^\circ=N_G(S)$, as needed.

Next, we have $p=11$. By \cite[Table I]{wilson1998}, $N_G(S)/S\cong 5\times 2A_5$, and $S\cong 11^2$. The centralizer of an element of order $11$ is $11\times M_{12}$. Thus for each subgroup $Q$ of order $11$, $K_G^\circ$ contains an element of order $5$ that centralizes $Q$ and normalizes a complement in $S$. Each of these twelve subgroups $5$ must be different, so $K_G^\circ$ contains all of $5\times 2A_5$. (Alternatively one may use that $K_G^\circ$ now contains $5\times 5$ and is normal in $N_G(S)$.)

Next, we have $p=7$, where $|N_G(S)/S|=36$. From \cite[p.234]{atlas}, the centralizers of 7A and 7B elements are $7\times He$ and $7^{1+4}.2A_7$, so we may apply Lemma \ref{lem:autocent}. If $N_G(S)$ is generated by $S$ and elements $x$ such that $|C_S(x)|>1$, then $K_G^\circ=N_G(S)$.

Since $N_G(S)$ is contained in a maximal subgroup $7^{1+4}_+\rtimes(3\times 2\cdot S_7)$ (this is the normalizer of a 7B element \cite[p.234]{atlas}, and contains $N_G(S)$ by \cite[Table II]{wilson1998}), we can easily calculate $C_S(x)$ for $x\in N_G(S)$. Of the $36$ elements $x$ in a complement $E$ to $S$ in $N_G(S)$, $24$ satisfy $|C_S(x)|>1$, and hence lie in $K_G^\circ$. Thus $K_G^\circ=N_G(S)$, as claimed.

For $p=5$, by \cite[Table I]{wilson1998}, $N_G(S)/S\cong 4\times 4\times S_3$, and it is contained in the normalizer of a 5B element by \cite[Table III]{wilson1998}. The centralizers of 5A and 5B elements are $5\times HN$ and $5^{1+6}_+.2J_2$ by \cite[p.234]{atlas}. In particular, we may apply Lemma \ref{lem:autocent}, as in the case $p=7$ above. This time there are again only twelve elements $x$ in a complement $E$ to $S$ in $N_G(S)$ such that $|C_S(x)|=1$, so again $K_G^\circ=N_G(S)$.

Finally, we have $p=3$. By \cite[Table I]{wilson1998}, $N_G(S)/S\cong 2\times 2\times SD_{16}$. As with the previous two cases, we apply Lemma \ref{lem:autocent}. The centralizers of 3A, 3B and 3C are $3Fi_{24}'$, $3^{1+12}_+.2\cdot Suz$ and $3\times Th$ respectively, so again we need to check that a complement $E$ to $S$ in $N_G(S)$ is generated by elements that are not fixed-point-free. We use the group $3^{2+5+10}.(M_{11}\times 2S_4)$ for this, as in \cite{wilson1998}, which contains $N_G(S)$. This time we find no fixed-point-free elements at all, and this completes the proof.

\subsection{Central extensions}
\label{sec:central}
We must also prove Theorem \ref{thm:sporadic} for central extensions of sporadic groups. In all cases but $2\cdot BM$, the central extension $\hat G$ can be constructed in Magma and we can test whether $Z(\hat G)$ is contained in some $O^{p'}(N_{\hat G}(Q))$, thus proving $Z(\hat G)\leq K_{\hat G}^\circ$ directly. If this is not the case, we may still have that $Z(\hat G)\leq K_{\hat G}^\circ$, and can easily test that on a computer as well. In the cases where $Z(\hat G)\not\leq K_{\hat G}^\circ$, we can also check whether $K_{\hat G}^\circ=K_{\hat G}$ using the same reasoning as with the simple case. For the various $\hat G$ and $p$, we give a choice of $Q$ for which $Z(\hat G)\leq O^{p'}(N_{\hat G}(Q))$ (when this is the case).

\subsubsection{All groups for the prime $2$}

We begin with $p=2$, so $|Z(\hat G)|=3$. There are six sporadic groups with Schur multiplier a multiple of $3$.
\begin{itemize}
\item If $\hat G=3\cdot M_{22}$, $Q$ may be taken to be a subgroup $2^2$ with normalizer of order $2^5\cdot 3^3$.
\item If $\hat G=3\cdot J_3$ there is no such $Q$, so $Z(\hat G)\not\leq K_{\hat G}^\circ$. However, $Z(\hat G)\leq K_{\hat G}$: this is because $z\in Z(\hat G)\setminus\{1\}$ may be expressed as a product $x_1x_2$ where $x_1\in O^{2'}(N_{\hat G}(Q_1))$ and $x_2\in O^{2'}(N_{\hat G}(Q_2))$. For $Q_1$ we choose a Klein four group with normalizer of order $2^7\cdot 3^3$ (there are two classes, and the other has normalizer of order $2^5\cdot 3^2$). For $Q_2$ we choose an elementary abelian group of order $16$ (unique up to conjugacy in $G$). We may choose $Q_1$ and $Q_2$ so that $Q_1\leq Q_2$, so any element of $N_G(S)$ that may be written as $x_1x_2$ lies in $K_{\hat G}$. The central $3$ is such an element. In the next section we confirm by direct computation that $T(\hat G,S)=1$.
\item If $\hat G=3\cdot McL$, $Q$ may be taken to be a subgroup $2^2$ with normalizer of order $2^7\cdot 3^3$.
\item If $\hat G=3\cdot Suz$, $Q$ may be taken to be the subgroup generated by a 2B involution with normalizer $2^{1+6}\cdot \PSU_4(2).2$ in the quotient $G$.
\item If $\hat G=3\cdot ON$, $Q$ may be taken to be the subgroup generated by an involution (all are conjugate), which has normalizer $4_2\cdot \PSL_3(4).2_1$ in the quotient $G$.
\item If $\hat G=3\cdot Fi_{24}'$, $Q$ may be taken to be the subgroup generated by a 2B involution (a $2$-central involution), which has normalizer $2^{1+12}_+\cdot 3\cdot \PSU_4(3).2$ in the quotient $G$.
\end{itemize}

\subsubsection{Odd primes and all groups except $2\cdot BM$}

If $p$ is odd then it is generally much easier to find a candidate subgroup $Q$ (as there are far fewer options for $Q$), and we simply list those groups $\hat G$ and primes $p$ such that $Z(\hat G)$ is contained in $O^{p'}(N_{\hat G}(Q))$.
\begin{itemize}
\item For $p=3$, if $\hat G$ is one of $2\cdot M_{12}$, $2\cdot J_2$, $2\cdot Suz$, $2\cdot Co_1$ and  $2\cdot Fi_{22}$, then there exists a non-trivial $p$-subgroup $Q$ such that $Z(\hat G)$ is contained in $O^{p'}(N_{\hat G}(Q))$.

If $\hat G=2\cdot HS$, then $Z(\hat G)$ is not contained in $K_{\hat G}^\circ$, and not in $K_{\hat G}$ by Lemma \ref{lem:nonon-cyclic}. However, in this case, the image of $Z(\hat G)$ in $N_{\hat G}(S)/K_{\hat G}$ is contained in the derived subgroup, so $T(\hat G,S)\cong T(G,S)$.

If $\hat G=4\cdot M_{22}$, then the central involution lies in $K_{\hat G}^\circ$, but the whole of the centre does not (hence not in $K_{\hat G}$ by Lemma \ref{lem:nonon-cyclic}). If $\hat G=2\cdot Ru$ then $Z(\hat G)$ does not lie in $K_{\hat G}^\circ$, but the proofs that these central elements also do not lie in $K_{\hat G}$ follow from the direct construction of elements of $T(\hat G,S)$ below.

\item For $p=5$, if $\hat G$ is one of $2\cdot J_2$, $2\cdot HS$, $2\cdot Suz$, $3\cdot Suz$ (and therefore $6\cdot Suz$), $2\cdot Ru$ and $2\cdot Co_1$, then there exists a non-trivial $p$-subgroup $Q$ such that $Z(\hat G)$ is contained in $O^{p'}(N_{\hat G}(Q))$.

For $3\cdot Fi_{24}'$, there is no subgroup $Q$ such that $Z(\hat G)\leq O^{5'}(N_{\hat G}(Q))$, but we still have that $Z(\hat G)\leq K_{\hat G}^\circ$.

For $2\cdot Fi_{22}$ and $3\cdot Fi_{22}$, and therefore $6\cdot Fi_{22}$, $Z(\hat G)$ does not lie in $K_{\hat G}^\circ$, which is equal to $K_{\hat G}$ by Lemma \ref{lem:nonon-cyclic}. The quotient $N_{\hat G}(S)/K_{\hat G}$ is $S_3\times Z(\hat G)$. In particular, $T(\hat G,S)\cong 2 \times Z(\hat G)$.

For $3\cdot McL$, since $S$ is TI, we have that $K_{\hat G}=S$, so $N_{\hat G}(S)/K_{\hat G}\cong N_G(S)/K_G\times Z(\hat G)$. In particular, $T(\hat G,S)\cong C_{24}$.

\item For $p=7$, if $\hat G$ is one of $2\cdot Co_1$ and $3\cdot Fi_{24}'$ then there exists a non-trivial $p$-subgroup $Q$ such that $Z(\hat G)$ is contained in $O^{p'}(N_{\hat G}(Q))$.

For $3\cdot ON$, there is no subgroup $Q$ such that $Z(\hat G)\leq O^{7'}(N_{\hat G}(Q))$, but we still have that $Z(\hat G)\leq K_{\hat G}^\circ$.

\end{itemize}

\subsubsection{The group $2\cdot BM$}

The only group with a central extension that isn't easy to work in is $G=BM$. The group $BM$ has a $4370$-dimensional representation over $\F_2$, and the dimensions of the fixed spaces of elements from 2A, 2B, 2C and 2D are $2510$, $2322$, $2212$ and $2202$ respectively. (This can be found by powering appropriate elements of large order in $BM$, which can be obtained using the words in the online Atlas \cite{webatlas}.) From \cite{atlas}, we see that (only) the 2C involution has preimage of order $4$ in $2\cdot G$, so if $O^{p'}(N_G(Q))$ possesses a 2C involution then the preimage in $2\cdot G$ is a non-split central extension.

The subgroup $5\times HS$ contains a 2C involution, so this deals with the case $p=5$. For $p=3$, the centralizer of a 6A element in $M$ is the intersection of the centralizer of a 2A element -- which is $2\cdot BM$ -- and a 3A element -- which is $3\cdot Fi_{24}'$. There are only two classes of involutions in $3\cdot Fi_{24}'$, and by checking orders we see that inside $Fi_{24}'$ the centralizer is $2\cdot Fi_{22}.2$. Thus inside $BM$, the preimage of $Fi_{22}$ in $2\cdot BM$ is $2\cdot Fi_{22}$, so there is a subgroup $3\times 2\cdot Fi_{22}$ in $2\cdot BM$, proving that $Z(\hat G)$ is in $K_{\hat G}^\circ$, as needed.

The centralizer of the $7$ is $7\times 2\cdot \PSL_3(4).2$, which is contained in the centralizer of an involution $2\cdot {}^2\!E_6(2).2$. We claim that the preimage of this centralizer of an involution in $2\cdot BM$ is $2^2\cdot {}^2\!E_6(2).2$. This, however, is easy to see: the subgroup $2^2\cdot {}^2\!E_6(2).S_3$ is the normalizer of a pure 2A subgroup of order $4$ in the Monster $M$, and therefore the centralizer of a point on the 2A subgroup in $M$ is $2^2\cdot {}^2\!E_6(2).2$. As this centralizes a 2A element, this lies inside the centralizer of a 2A element of $M$, namely $2\cdot BM$.

As there is a copy of $S_3$ permuting the central involutions of the $2^2\cdot {}^2\!E_6(2)$, all central extensions $2\cdot {}^2\!E_6(2)$ are isomorphic. As the centralizer of an element of order $7$ in $2\cdot {}^2\!E_6(2).2$ is $7\times 2\cdot \PSL_3(4).2$, the preimage in $2^2\cdot {}^2\!E_6(2).2$, hence in $2\cdot BM$, of this subgroup must be $7\times 2^2\cdot \PSL_3(4).2$, as otherwise not all of the quotients by a central involution would be isomorphic. In particular, this means that $Z(\hat G)$ lies inside $O^{7'}(N_{\hat G}(Q))$, and therefore $K_{\hat G}^\circ$ contains $Z(\hat G)$, as needed.

\subsection{Outer automorphisms}
\label{sec:autos}

This section details the actions of the outer automorphism of the simple group $G$ (if it has one) on $N_G(S)/K_G$, thereby computing the group $N_{\Aut(G)}(S)/K_G$. We of course need only do this when $K_G\neq N_G(S)$, and when $G\neq \Aut(G)$. We include the cases where $S$ is TI, even though then it is simply $N_{\Aut(G)}(S)/S$.

\begin{enumerate}
\item For $M_{22}$ and $p=3$, $N_{\Aut(G)}(S)/K_G\cong SD_{16}$.
\item For $J_2$, if $p=3$ then $N_{\Aut(G)}(S)/K_G\cong 2\times 2$, and if $p=5$ then $N_{\Aut(G)}(S)/K\cong 4$.
\item For $J_3$, if $p=3$ then $N_{\Aut(G)}(S)/K_G\cong 2\times 2$.
\item For $HS$, if $p=3,5$ then $N_{\Aut(G)}(S)/K_G\cong D_8$.
\item For $McL$, if $p=5$ then $N_{\Aut(G)}(S)/K_G\cong 3\rtimes E$ for a $2$-group $E$, with $C_4\times C_2$ as quotient. (The group $E$ has order $16$ and a cyclic subgroup of index $2$, but is not of maximal class. This determines $E$ uniquely.)
\item For $Suz$, if $p=5$ then $N_{\Aut(G)}(S)/K_G\cong 2\times 2$.
\item For $He$, if $p=5$ then $N_{\Aut(G)}(S)/K_G\cong S_3$.
\item For $Fi_{22}$, if $p=5$ then $N_{\Aut(G)}(S)/K_G\cong D_{12}$.
\item For $T.2$, if $p=5$ then $N_{\Aut(G)}(S)/K_G\cong 4S_4$, which has $C_4$ as quotient.
\end{enumerate}

\section{Constructing endotrivial modules}

We attempt to identify the trivial-source endotrivial modules where there are non-trivial ones, and the simple endotrivial modules that were not determined in \cite{lassueurmalleschulte2016}.

The characters of many of these have previously been identified in \cite{lassueurmazza2015}, but here we give the full socle series of the modules. The simple trivial-source, endotrivial modules that are given in Section \ref{sec:trivialsource} were all found already in \cite{lassueurmalleschulte2016}.

\subsection{Trivial-source endotrivial modules}
\label{sec:trivialsource}
For the small sporadic groups, we simply induce up the modules for $N_G(S)/K_G$ to $G$, passing through a maximal subgroup between $N_G(S)$ and $G$ when one is available. We then remove all summands not from the appropriate block (the Brauer correspondent of the block containing the module for $N_G(S)$), and then remove any projective summands if we can find them. At this point, we just ask Magma to write the module as a direct sum.

For larger modules, we are sometimes able to write the induction (cut by the appropriate block) as a sum of a projective and another module: if Magma can construct the projective module $P$, the command \texttt{AHom(P,V)} constructs all $kG$-module homomorphisms from $P$ to $V$, so we may take the cokernel of an injective map in this space. In such cases where we wish to do this, we give a direct sum decomposition of the induction to avoid the reader having the rediscover the projective part of the induction.

\subsubsection{The groups $M_{11}$, $M_{22}$ and $M_{23}$}
In all three cases we have $p=3$. For $G=M_{11}$ we have $T(G,S)\cong C_2\times C_2$. The four modules are $1$, $10_1$ (the self-dual $10$-dimensional module), a module
\[ 10_2/5/1,24/5^*/10_2^*,\]
and a module
\[ 5^*/((5/1,24/5^*)\oplus 10_1)/5.\]

For $G=M_{22}$ and $\hat G=2\cdot M_{22}$ we have $T(G,S)\cong C_2\times C_2$ and $T(\hat G,S)\cong C_4\times C_2$. The four modules in $T(G,S)$ are $1$, $55$,
\[ 49/1,55/49^*,\quad 49^*/1,55/49.\]
The four faithful modules in $T(\hat G,S)$ are all simple, and are $10$, $10^*$, $154$ and $154^*$.

For $G=M_{23}$ we have $T(G,S)\cong C_2$, and the two modules are $1$ and $253$.

\subsubsection{The groups $J_2$, $J_3$ and $J_4$}

For $G=J_2$ we have $p=3,5$, and in both cases $T(G,S)\cong C_2$. The non-trivial module in $T(G,S)$ for $p=3$ and $p=5$ have structures
\[ 21_1,21_2/133/21_1,21_2\quad\text{and}\quad 189/14,21,189/14,21,189/189\]
respectively.

For $G=J_3$ and $\hat G=3\cdot J_3$, if $p=3$ then $T(G,S)\cong C_2$. The two modules are $1$ and a summand $V$ of the $25840$-dimensional induced module from $N_G(S)$. This induced module is the sum of the projectives $P(2754)=2754/324/2754$, $1215_1$ and $1215_2$, the simple module $324$ and $V$. The module $V$ has dimension $17254$, $85$ composition factors, $25$ socle layers, and socle structure
\[ 168/306,934/36,168/306/168/306/36,168/1^2,306/168/306,934/36,168/\]
\[1^2,306/36,168/1^4,306,934^2/36,168^2/1^4,306^2/168^2/1,306^2,934/36^2,168^3/1^4,306^2/\]
\[36^2,168^3/1^4,306^2,934^2/1^3,36,168^2,934/1^4,36,168,306/168,934.\]
(Here, $934^2$ means two copies of $934$, and so on.)

If $p=2$ then $T(\hat G,S)$ is trivial. Since $K_G\neq K_G^\circ$, and \cite{lassueurmazza2015} is in disagreement, we confirm directly that $K_{\hat G}=N_{\hat G}(S)$. To do this we induce one of the $1$-dimensional characters of the group $3\times (2^5\cdot A_5)$ up to all of $\hat G$, as suggested in \cite[4.2]{lassueurmazza2015}. The summands have dimensions $12681$, $4770$ and $8712$, so the one of dimension $12681$ is the Green correspondent. However, its restriction to $S$ is not the sum of a trivial and a free module, as would be needed. Indeed, upon removing the free summands, there is a module of dimension $1161$ left over. The problem in \cite{lassueurmazza2015} appears to have been a communication problem between Carlson and the authors of \cite{lassueurmazza2015}, and Carlson's algorithm produces the correct answer.

For $J_4$ and $p=11$, $N_G(S)$ is maximal, and so it is very difficult to obtain any information about the elements of $T(G,S)$, beyond the fact that it is isomorphic to $C_{10}$. In \cite{lassueurmalleschulte2016}, candidate simple endotrivial modules are given as the two dual modules of dimension $887778$, and the self-dual module of dimension $394765284$. In the next section we show that the latter cannot be endotrivial, and that the former is endotrivial, but does not have trivial source. We cannot construct any non-trivial element of $T(G,S)$.

\subsubsection{The group $HS$}
Let $G=HS$. If $p=3$, we have $T(G,S)\cong C_2\times C_2$. The four modules are $1$, $154_1$, $154_2$ and $154_3$.

If $p=5$, we have $T(G,S)\cong C_4$. Choose one of the two dual generators $V_1$ for $T(G,S)$, and let $V_2$ denote the element of order $2$ from $T(G,S)$. Let $W_i$ denote the induction of the Green correspondent of $V_i$ in $N_G(S)$ to $G$. Thus $\dim(W_i)=22176$. We may decompose $W_1$ as
\[ V_1\oplus P(280^*)\oplus P(650)\oplus P(1925)\oplus 1925 \oplus 1750.\]
The module $V_1$ has dimension $2876$, and has socle series
\[98/1,133_1,133_2/21,98/55,55,210,518/21,21,98,98,98,280/1,55,133_1,133_2,518/98.\]
(The heart of this module is indecomposable, and a Hasse diagram for it is too complicated to draw nicely.) Its dual has the same structure except for $280^*$ in place of $280$.

The other module is a self-dual $12376$-dimensional module $V_2$ with $64$ composition factors and seven socle layers. The decomposition of $W_2$ is
\[ V_2\oplus P(1925)\oplus 175\oplus 1750\oplus 2750,\]
and the precise socle structure of $V_2$ is
\[
210,518^2/21^2,98^2,280_1,280_2/
55,133_1,133_2,210^2,518^2/
21^4,98^5,280^2_1,280^2_2/\]
\[1^4,55^2,133_1^2,133^2_2,210^3,518^5/21^4,98^5,280^2_1,280^2_2/55,210,518^2.\]

\subsubsection{The group $McL$}
\label{subsubsec:mcl}
Let $G=McL$ and $p=5$. Writing $\hat G=3\cdot McL$, we have that $T(G,S)\cong C_8$ and $T(\hat G,S)\cong C_{24}$. There are four simple endotrivial modules for $\hat G$, all of dimension $126$. However, their sources are $\Omega^{\pm 2}(k)$, so they do not lie in $T(\hat G,S)$. On the other hand, the tensor product of the two $126$-dimensional modules in the same block is indecomposable and has trivial source, so does lie in $T(\hat G,S)$. It has socle structure
\[1035/126_2/1233/126_1/1035/126_2/639,1233/126_1,126_1/1035,4752/126_2,126_2,153,846/639,1233/126_1/1035.\]
By finding its Green correspondent in $N_{\hat G}(S)$, we see that it has order $6$ in $T(\hat G,S)$. It seems difficult to construct the other non-trivial elements of $T(\hat G,S)$.

\subsubsection{The groups $Suz$, $He$, $Co_3$, $Co_2$ and $Fi_{23}$}
In these cases $p=5$. If $G=Suz$ then $T(G,S)\cong C_2$, and the two modules are $1$ and $1001$. If $G=He$ then  $T(G,S)\cong C_3$, and the three modules are $1$, $51$ and $51^*$. If $G=Fi_{23}$, then $T(G,S)\cong C_2$, and the two modules are $1$ and $111826$.

If $G=Co_3$ then $T(G,S)\cong C_4$. We cannot construct the two elements of order $4$, but we can construct the element of order $2$, by inducing the non-trivial $1$-dimensional module from the maximal subgroup $McL\rtimes 2$. (Note that the other trivial-source endotrivial modules can be induced from modules for $McL\rtimes 2$, but we cannot construct them either: see Section \ref{subsubsec:mcl}.) The induction has structure $23/230/23$.

If $G=Co_2$ then $T(G,S)\cong C_2$. Unfortunately we cannot construct the non-trivial element of $T(G,S)$ in this case, but it also comes from a non-trivial $1$-dimensional module, this time for the maximal subgroup $HS\rtimes 2$. This induced module $W$ has dimension $476928$, and from an ordinary character calculation, and the known decomposition matrix for $Co_2$, we see that the contribution $V$ from the principal $5$-block has composition factors $23,23,2254,2254,29624$. By restricting $23$ and $2254$ down to $HS\rtimes 2$ and using the Nakayama relation, we find that $\soc(V)=23\oplus 2254$. However, this is not enough to give us the complete structure of $V$, since $\Ext^1(23,2254)\neq 0$, so $V$ need not have structure $23,2254/29624/23,2254$.

\subsubsection{The group $Ru$}
For $G=Ru$ and $\hat G=2\cdot Ru$, we have $p=3$. We have $T(G,S)\cong C_2$ and $T(\hat G,S)\cong C_4$. The elements of order $4$ in $T(\hat G,S)$ are the modules $28$ and $28^*$; the symmetric square of $28^\pm$ is $406$, which is the non-trivial module in $T(G,S)$.

\subsubsection{The group $Fi_{22}$}
\label{sec:fi22const}

For $Fi_{22}$, $p=5$, we have $T(G,S)\cong C_2$. The two modules are $1$ and $1001$.

In general, $T(\hat G,S)\cong C_2\times Z(\hat G)$. If $|Z(\hat G)|$ has order $2$ there are two faithful trivial-source endotrivial modules. Both have dimension $61776=|G:H|$ for $H$ the maximal subgroup $Z(\hat G)\times \Omega_8^+(2).S_3$, and are the inductions of $1$-dimensional modules from this maximal subgroup. The two modules have structure
\[13376_i/(352\oplus(5824_i/23024_i/5824_i^*))/13376_i\]
for $i=1,2$, with an appropriate labelling of the simple modules. (Note that this contains $5824_i^\pm$, which we need in the next subsection.)

For $\hat G=3\cdot Fi_{22}$, the two $351$-dimensional simple modules are endotrivial and trivial source, being the elements of order $3$ in $T(\hat G,S)$. The two elements of order $6$ are dual and also must be Galois conjugates, so are `Galois dual'. The composition factors, from an ordinary character calculation and the decomposition matrix on the online Modular Atlas (determined from \cite{hisswhite1994}), are
$2079, 5643^2, 48411$. By Galois duality, the module must be $5643/2079,48411/5643$.

In the case where $Z(\hat G)$ has order $6$, the four $61776$-dimensional modules are simple and endotrivial. This is easy to check because the inductions of the relevant $1$-dimensional modules for $H$ are, in fact, simple.

\subsubsection{The groups $Th$ and $BM$}

For $Th$ and $BM$, $p=7$, there are two elements in $T(G,S)$, but it is impossible with our current technology to construct the non-trivial element of it. (However, it is easy to understand the corresponding linear character of $N_G(S)$, i.e., the Green correspondent of that element of $T(G,S)$.)

Note that if $G=BM$ then there is a subgroup $H\leq G$ isomorphic with $Th$. Indeed, up to conjugacy $N_G(S)\leq H$, and therefore the Green correspondents of the elements of $T(H,S)$ are the elements of $T(G,S)$.

Furthermore, since the maximal subgroup $L\cong 2^2\cdot F_4(2).2$ contains $N_G(S)$, we see that the Green correspondents of the elements of $T(G,S)$ are the two linear characters of $L$.


\subsubsection{The group $T$}

For the Tits group $T$, $p=5$, we have $T(G,S)\cong C_6$. A generator for this group is the simple module $26$ (or $26^*$) with its tensor square splitting as the sum of the simple module $325$ and a $351$-dimensional module in $T(G,S)$, which has structure
\[ 109_1/27/1,78/27^*/109_2,\]
for some suitable labelling of $109_i$ and $27$. 

For the final module we may work over $\F_5$, so that $109_1\oplus 109_2$ becomes $218$ and $460_1\oplus 460_2$ becomes $920$. The induction $W$ from $N_G(S)$ to $G$ has dimension $14976$, and splits up as
\[ V\oplus P(920)\oplus 300^{\oplus 2}\oplus 1300_1\oplus 1300_2,\]
where $V$ is the trivial-source endotrivial module, of dimension $2976$. It has socle structure
\[27/920/27^*/1,351/27^*,52/218,351^*/27,27/1,920/27^*.\]
(Here we have chosen $27$ so that there is a module $1/27^*$ but not $1/27$, and $351$ so that $\Lambda^2(27)=351^*$.)

\subsection{Simple endotrivial modules}

This section tabulates the known simple endotrivial modules for sporadic groups, and gives their sources. It also gives the candidate simple endotrivial modules from \cite{lassueurmalleschulte2016} whose status we have not determined here.

At the time of writing, the known simple endotrivial modules are given in Table \ref{tab:knownsimple} (if the module was known before this paper, a reference is given in the table), and the remaining candidates for simple endotrivial modules are in Table \ref{tab:candidatesimple}.

\begin{table}\begin{center}\begin{tabular}{ccccc}
\hline Group & Prime & Module & Source & Previously known?
\\\hline $M_{11}$ & $3$ & $10_1$ & $k$ & \cite{lassueurmalleschulte2016}
\\ $M_{11}$ & $3$ & $10_2^\pm$ & $\Omega^{\pm 2}(k)$ & \cite{lassueurmalleschulte2016}
\\ $M_{22}$ & $3$ & $55$ & $k$ & \cite{lassueurmalleschulte2016}
\\ $2\cdot M_{22}$ & $3$ & $10^\pm$ and $154^\pm$ & $k$ & \cite{lassueurmalleschulte2016}
\\ $M_{23}$ & $3$ & $253$ & $k$ & \cite{lassueurmalleschulte2016}
\\ $HS$ & $3$ & $154_1$, $154_2$ and $154_3$ & $k$ & \cite{lassueurmalleschulte2016}
\\ $3\cdot McL$ & $5$ & $126_1^\pm$ and $126_2^\pm$ & $\Omega^{\pm 2}(k)$ & \cite{lassueurmalleschulte2016}
\\ $Suz$ & $5$ & $1001$ & $k$ & \cite{lassueurmalleschulte2016}
\\ $He$ & $5$ & $51^\pm$ & $k$ & \cite{lassueurmalleschulte2016}
\\ $Ru$ & $3$ & $406$ & $k$ & \cite{lassueurmalleschulte2016}
\\ $2\cdot Ru$ & $3$ & $28^\pm$ & $k$ & \cite{lassueurmalleschulte2016}
\\ $3\cdot ON$ & $7$ & $342_1^\pm$ and $342_2^\pm$ & $U_1$ & \cite{lassueurmalleschulte2016}
\\ $Fi_{22}$ & $5$ & $1001$ & $k$ & \cite{lassueurmalleschulte2016}
\\ $2\cdot Fi_{22}$ & $5$ & $5824_1^\pm$ and $5824_2^\pm$ & $\Omega^{\pm 5}(k)$ & No
\\ $3\cdot Fi_{22}$ & $5$ & $351^\pm$ & $k$ & \cite{lassueurmalleschulte2016}
\\ $3\cdot Fi_{22}$ & $5$ & $12474_1^\pm$ and $12474_2^\pm$ & $\Omega^{\pm 5}(k)$ & No
\\ $6\cdot Fi_{22}$ & $5$ & $61776_1^\pm$ and $61776_2^\pm$ & $k$ & \cite{lassueurmalle2015}
\\ $Fi_{23}$ & $5$ & $111826$ & $k$ & \cite{lassueurmalleschulte2016}
\\ $J_4$ & $11$ & $887778$ & $\Omega^{\pm 8}(k)$ & No
\\ $T$ & $3$ & $26^\pm$ & $U_2$ & \cite{lassueurmalleschulte2016}
\\ $T$ & $5$ & $26^\pm$ & $k$ & \cite{lassueurmalleschulte2016}
\\ $T$ & $5$ & $351^\pm$ & $\Omega^{\pm 6}(k)$ & \cite{lassueurmalleschulte2016}
\\ \hline
\end{tabular}\end{center}
\caption{Known simple endotrivial modules for quasisimple sporadic groups. For definitions of $U_1$ and $U_2$, see the proof of Proposition \ref{prop:knownsimple}.}\label{tab:knownsimple} 
\end{table}

\begin{table}\begin{center}\begin{tabular}{cccc}
\hline Group & Prime & Module & Notes
\\\hline $Th$ & $7$ & $27000^\pm$ & Known to be simple
\\ $BM$ & $7$ & $9287037474$ & Simplicity unknown, trivial source if endotrivial
\\ $BM$ & $7$ & $775438738408125$ & Simplicity unknown, trivial source if endotrivial
\\ \hline
\end{tabular}\end{center}
\caption{Remaining candidate simple endotrivial modules for quasisimple sporadic groups.}\label{tab:candidatesimple} 
\end{table}

\begin{prop}\label{prop:knownsimple} If $V$ appears in Table \ref{tab:knownsimple}, then $V$ is a simple endotrivial module.
\end{prop}
\begin{pf} If $V$ appears in \cite[Table 7]{lassueurmalleschulte2016} without a `?' then $V$ was shown to be endotrivial in that paper. The sources of simple modules of small dimension for small $G$ can easily be found by computer. If $V$ is self-dual then it must be trivial source as $p$ is odd, and thus $T(S)$ is torsion-free by \cite{carlsonthevenaz2004}, as we stated in Section \ref{sec:prelim}.

Of those simple modules whose endotriviality was established in \cite{lassueurmalleschulte2016}, the only ones that are difficult to construct (because they do not have very small dimension and are unavailable in the online Atlas \cite{webatlas}) are the $342$-dimensional modules for $3\cdot ON$. One may find these simple modules as composition factors of a tensor product $45\otimes 495_i$, and so their source can be accessed relatively easily. It can be described as follows: there are three conjugacy classes of subgroups $Q_i$ of order $49$ in $N_G(S)$, with (say) $Q_1$ and $Q_2$ being exchanged by the outer automorphism. The action of both $Q_1$ and $Q_3$ is (up to projectives) $\Omega^3(k)$, and on $Q_2$ it is $\Omega^{-11}(k)$. This determines an endotrivial module $U_1$ of $S$ uniquely up to outer automorphism of $G$, by Lemma \ref{lem:restrict}.

For the Tits group, if $p=5$ then the sources are as in Table \ref{tab:knownsimple}, but for $p=3$ the source is more difficult to describe, as it was for the O'Nan group. Let $Q_1$ and $Q_2$ be representatives of the two $N_G(Q)$-classes of subgroups $3^2$ in $S$. The module $U_2$ has the property that $\Omega^3(U_2)$ has restriction to $Q_1$ trivial plus projective, and $\Omega^{-3}(U_2)$ has restriction to $Q_2$ trivial plus projective. This determines $U_2$ uniquely up to outer automorphism of $G$, which swaps $Q_1$ and $Q_2$.

We are left with those modules that were not proved to be endotrivial in \cite{lassueurmalleschulte2016}: these are the modules for $2\cdot Fi_{22}$ and $6\cdot Fi_{22}$. The $61776$-dimensional modules $V$ are, by an ordinary character calculation, induced from $1$-dimensional modules for the maximal subgroup $H\cong 6\times \Omega_8^+(2).S_3$, which is of course of index $61776$. Thus these are trivial-source modules. For endotriviality we either note that they arise from characters on $N_G(S)/K_G$, or we note that by the Mackey formula we need to understand $H^g\cap S$ for $g$ running over all 2472 $(H,S)$-double-coset representatives. All but one of these intersections is trivial, so $V{\downarrow_S}$ is the sum of a trivial and a free module.

For $V$ one of $5824_i^\pm$ for $\hat G=2\cdot Fi_{22}$, one may find these simple modules in the elements of $T(\hat G,S)$ that were given in the previous section. Once they are given in a computer it is easy to restrict to $S$ and remove free summands. One is left with $\Omega^{\pm 5}(k)$, so $V$ is endotrivial with $\Omega^{\pm 5}(k)$ as source.

For $V$ one of the $12474^\pm_i$ for $\hat G=3\cdot Fi_{22}$, it is much harder to construct this module, and I am grateful to Richard Parker for running this computation on Meataxe64. One finds a $12474$-dimensional composition factor in the tensor product $5864\otimes 78$ over $\F_{25}$. Restricting down to $S$ and removing free summands, we are again left with $\Omega^{\pm 5}(k)$.

For $G=J_4$, we can confirm directly that the exterior squares of the $1333$-dimensional modules are endotrivial. This can be proved by restricting to the (unique up to conjugacy) maximal elementary abelian subgroup of order $11^2$. One cannot restrict the $887778$-dimensional module directly, so we restrict the $1333$-dimensional module, remove any free summands -- which reduces the dimension to $244$ -- then take the exterior square. This is the direct sum of a free module and the module $\Omega^{\pm 8}(k)$, so the exterior square of the $1333$-dimensional module for $G$ is endotrivial by Lemma \ref{lem:restrict}. I am grateful to Thomas Breuer who proved simplicity of the reduction modulo $11$ of this ordinary character: using GAP \cite{GAP4}, he showed that there is no consistent breaking up of the mod $11$ reduction of the ordinary character to two maximal subgroups $2^{11}\rtimes M_{24}$ and $2^{1+12}.3.M_{22}\rtimes 2$.
\end{pf}

\begin{prop} Let $\hat G$ be a quasisimple sporadic group and $V$ is a non-trivial, simple $k\hat G$-module. If $V$ does not appear in Tables \ref{tab:knownsimple} or \ref{tab:candidatesimple} then $V$ is not endotrivial.
\end{prop}
\begin{pf} We start with \cite[Table 7]{lassueurmalleschulte2016}, and therefore need to eliminate the modules that appear in that table but not in Tables \ref{tab:knownsimple} and \ref{tab:candidatesimple}.

First, notice that if $V$ is self-dual and endotrivial then $V$ lies in $T(\hat G,S)$. Thus if $T(\hat G,S)=1$, or is non-trivial but $V$ does not appear in the constructions in the previous section then $V$ is not endotrivial. This eliminates the module of dimension $74887473024$ for $Fi_{24}'$ and the module of dimension $7226910362631220625000$ for $M$, but does not necessarily eliminate the module of dimension $394765284$ for $J_4$ for $p=11$. (The modules for $M$ and $Fi_{24}'$ were also eliminated in \cite{lassueurmalle2015}.)

To do this, since in this case $V$ is self-dual, it must be the unique non-trivial self-dual element of $T(G)\cong \Z\times C_{10}$. Assuming $V$ is endotrivial, the Green correspondent of $V$ is the non-trivial real linear character of $N_G(S)$, and this takes value $-1$ on one class of elements of order $22$ in $G$. However, from \cite{atlas}, we see that $V$ takes value $1$ on both classes of elements of order $22$ in $G$. Since projective characters take value $0$ on $p$-singular classes, this means that the Green correspondent of this linear character cannot be the reduction modulo $11$ of this ordinary character, and so the module $V$ of dimension $394765284$ cannot be endotrivial.
\end{pf}

We cannot do anything with the modules for $Th$, since they are not self-dual, and therefore they could just be elements of infinite order in $T(G)$. For $BM$, at most one of the two modules listed is endotrivial, because it would have to lie in $T(G,S)$ for this group, which has order $2$.


\bibliographystyle{amsplain}
\bibliography{references}
\end{document}